\documentclass[10pt,letterpaper]{article}
\usepackage{t1enc}
\usepackage[latin1]{inputenc}
\usepackage[english]{babel}
\usepackage{graphics}
\begin{document}
\date{}
\title{  Expansions of Theta Functions and Applications}
\author{A.  Raouf  Chouikha \footnote
{Universite Paris 13 LAGA UMR 7539 Villetaneuse 93430,  e-mail: 
chouikha@math.univ-paris13.fr}
}
\maketitle

\begin{abstract} We prove that the classical theta function $\theta_4$ 
may be expressed as 
$$ \theta_4(v,\tau) = \theta_4(0,\tau)\ \exp[- \sum_{p\geq 1} \sum_{ 
k\geq 0} \frac {1}{p} \bigg( \frac {\sin \pi v}{(\sin (k+\frac {1}{2})\pi 
\tau)}\bigg)^{2p}].$$ We obtain an analogous expansion for the three 
other theta functions since they are related. \\    
  These results have several consequences. In particular, an expansion 
of the Weierstrass elliptic function will be derived. Actions of the 
modular group and other arithmetical properties will also be considered. 
Finally using a new expression for the Rogers-Ramanujan continued 
fraction we produce a simple proof of a Rogers identity.

\end{abstract}

\section{Introduction}
The classical theta function is given by the doubly infinite sum
\begin{equation} \theta (v,\tau ) =  \sum_{-\infty < n < +\infty } 
e^{i\pi n^2\tau} e^{2i\pi nv}.
\end{equation}
which locally converges uniformly for $v \in C$ the complex plane, and 
$\tau \in {\cal H}_+$ the upper-half plane of complex numbers with 
positive imaginary part. \\
A remarquable feature of the theta function is its dual nature. Indeed, 
when viewed as a function of $v,$ we see it as an elliptic analog of 
the exponential function and it may be used in order to express elliptic 
functions. Since $\theta $ is periodic with period $1$ and quasi-period 
$\tau .$ When considered as a function of $\tau ,\ \theta $ has modular 
properties with close connection to the partition function and the 
representation of integers as sums of squares.\\
The story of this function started with Euler, J. Bernoulli and Fourier 
who occasionnaly used the theta function among  other closely related 
functions. However, the systematic study of theta functions and their 
utilisation for the theory of elliptic functions is due to Jacobi. He 
remarked that theta can be represented by series whose convergence is very 
fast and which may used for numerical computations of elliptic 
functions. \\  
General properties of theta functions may be found in many references, 
among which [B], [M] or [W-W] contain a good description. \\

The above function $\theta $ is often called the $\theta _3$ function 
of Jacobi. There are three other $\theta $ functions obtained by a 
change of characteristics
$$\theta _1(v,\tau ) =\sum_{-\infty \leq n\leq +\infty}(-1)^{\frac 
{n-1}{2}} q^{(\frac {n+1}{2})^2} e^{i\pi (2n+1)v} = 2 \sum_{n\geq 0} (-1)^n 
q^{(\frac{n+1}{2})^2} \sin ((2n+1) \pi v)$$
$$\theta _2(v,\tau ) =\sum_{-\infty \leq n\leq +\infty } q^{(\frac 
{n+1}{2})^2} e^{i\pi (2n+1)v} = 2 \sum_{n\geq 0} q^{(\frac{n+1}{2})^2} \cos 
((2n+1) \pi v)$$
$$\theta _4(v,\tau ) =\sum_{-\infty \leq n\leq +\infty } (-1)^n q^{n^2} 
e^{i\pi 2nv} = 1 + 2 \sum_{n\geq 1 } (-1)^n q^{n^2} \cos (2n \pi v)$$
where \ $q = e^{i \pi \tau }$\ verifies \ $\mid q\mid < 1,$\ and \ $v$\ 
is a complex number.\\
The four theta functions are related.
Then we can choose any of these four theta functions and then define 
the remaining three in terms of the one chosen.\\
All four theta functions are entire functions of $v$. All are periodic, 
the period of \ $\theta _1$\ and  \ $\theta _2$\ is\ $2$, and that of  
\ $\theta _3$\ and  \ $\theta _4$\ is\  $1$.\\
It is known that the zeros of $\theta_1$ are $m + n\tau $, those of 
$\theta_2$ are $\frac {1}{2} + m + n\tau $, those of $\theta_3$ are $\frac 
{1}{2} + \tau + m + n\tau $ and those of $\theta_4$ are $\frac 
{1}{2}\tau  + m + n\tau $ where $m,n$ are integers.\\
From the knowledge of the zeros it is possible to obtain infinite 
products representing the theta functions, and from these products the 
partial fraction expansions of $log\ \theta (v,\tau )$ follow.\\

\bigskip

The main goal of this paper is to state a trigonometric expansion of 
theta functions in powers of $\sin\pi v$ and to derive some applications. 
Namely $$ \log\frac{\theta_4(v,\tau)}{\theta_4(0,\tau)} 
=\sum_{p=1}^\infty c_{2p}(\tau)(\sin\pi v)^{2p}$$
where the coefficients $c_{2p},$ which depend on $\tau ,$ have the 
following form
$$c_{2p}(\tau ) = (-1)^{p+1} \frac {2^{2p+1}}{2p} \sum_{k \geq 0}\bigg( 
\frac {q^{2k+1}} {(1 - q^{2k+1})^{2}}\bigg)^p = - \frac {1}{p} 
\sum_{k\geq 0} \frac {1}{(\sin (k+\frac {1}{2})\pi \tau)^{2p}}$$ 
where \ $q = e^{i \pi \tau }, \quad \mid q\mid < 1.$\ 
 In the same way, we also obtain the expressions
$$\theta_3(v,\tau) = \theta_4(0,\tau)\ \exp[- \sum_{p\geq 1} \sum_{ 
k\geq 0}\frac {1}{p} \bigg( \frac {\cos \pi v}{(\sin (k+\frac {1}{2})\pi 
\tau)}\bigg)^{2p}],$$
$$\theta_2(v,\tau) = \theta_4(0,\tau)\ \exp[ i\pi (v+\frac {1}{4}\tau)- 
\sum_{p\geq 1} \sum_{ k\geq 0}\frac {1}{p} \bigg( \frac {\cos \pi 
(v+\frac {1}{2}\tau)}{(\sin (k+\frac {1}{2})\pi \tau)}\bigg)^{2p}],$$
$$\theta_1(v,\tau) = \theta_4(0,\tau)\ \exp[ i\pi (v-\frac {1}{2}+\frac 
{1}{4}\tau)- \sum_{p\geq 1} \sum_{k\geq 0}\frac {1}{p} \bigg( \frac 
{\sin \pi (v+\frac {1}{2}\tau)}{(\sin (k+\frac {1}{2})\pi 
\tau)}\bigg)^{2p}].$$ 
The above expressions for \ $\theta_4 $ \ and \ $ \theta_3$ \ are valid 
in the "strip" \ $\mid \frac {\sin \pi v}{\sin \frac {1}{2}\pi \tau} 
\mid < 1    ,$\\ those relating to  \ $\theta_2 $ \ and \ $ \theta_1$ \  
are valid in the "strip" \ $\mid \frac {\sin \pi (v+\frac {1}{2}\tau 
)}{\sin \frac {1}{2}\pi \tau} \mid < 1  .$. \\
 
Some consequences of the above expansions may be deduced. In 
particular, the Weierstrass elliptic function $\wp (z)$ with primitive periods \ 
$2, 2\tau $ \ has an analogous expansion
$$\wp(z+\tau) = e_3 - \sum_{p\geq 1}[  \sum_{k\geq 0} \frac {- 
2(2p+1)(\sin{\pi z\over {2}})^{2p}}{(\sin (k+\frac {1}{2})\pi \tau)^{2p+2}} +  
\sum_{k\geq 0} \frac {4p(\sin{\pi z\over {2}})^{2p}}{(\sin (k+\frac 
{1}{2})\pi \tau)^{2p}} ].$$
 We also examine properties of $c_{2p}$ under the actions of the 
modular group and Landen or Gauss transformations.  Some arithmetical 
applications using Lambert series  will also be considered.
An expression for the Rogers-Ramanujan continued fractions will be 
derived
$$ R(q) = q^{2/5} e^{ \sum_{p\geq 1,k\geq 0}{\displaystyle 
\frac{q^{2p(k+1)}}{p}\bigg(\frac {(q^2-q)^{2p}- 
(q^3-1)^{2p}}{(q^{5(2k+1)}-1)^{2p}}\bigg)}}.$$
The last one allows us to prove the well known Rogers identity
$$R(q) = q^{\frac {2}{5}} \prod_{k\geq 1} \frac {(1- q^{10k-2})(1- 
q^{10k-8})}{(1- q^{10k-4})(1- q^{10k-6})}.$$

\newpage

\section{Expansions of the theta functions}
\subsection {Classical expansions}
Recall at first the following facts which will be useful in the sequel
\bigskip
  
{\bf Proposition 2-1 }\quad {\it The function \ $\theta_4 (v,\tau )$  
satisfies the triple product
$$\theta_4 (v,\tau ) = \prod _{n\geq 1}(1 - q^{2n})(1-q^{2n+1}e^{2i\pi 
v})(1-q^{2n-1}e^{-2i\pi v})$$
 where $q = e^{i\pi \tau},\ v \in C,\ \tau \in {\cal H}_+.$}

\bigskip

It is also known that $log\ \theta_4$ has a Fourier series expansion, 
[B] or [W-W]
\bigskip

{\bf Proposition 2-2 }\quad {\it The theta function $\theta_4$ has the 
following expansion 
$$ \theta_4(v,\tau) = \theta_4(0,\tau)\ exp[4 \sum_{n³1} 
{q^n\over{1-q^{2n}}} {(\sin{n\pi v})^2\over{n}}].$$
Moreover, this expansion is valid for $Im \tau > 0$ and for $\mid Im 
v\mid < \frac {1}{2} Im \tau . $} 

	\bigskip 
We also obtain similar expansions for the other three theta functions, 
[B].\\  
In [C] we proved the following

\bigskip

{\bf Proposition 2-3 }\quad {\it The theta function $\theta_4$ may be 
expressed as 
$$\theta_4(v,\tau) = \theta_4(0,\tau)\ \exp \ [ \sum_{p\geq 1}
c_{2p}(\tau) (\sin\pi v)^{2p}]$$
  where the coefficients \ $c_{2p}$ \ satisfy
the recursion relation 
\begin{equation}
(2p+2)(2p+1) c_{2p+2}(\tau) - 4p^2 c_{2p}(\tau) =  \frac {4}{ \pi ^2} 
a_{2p}(\tau ).
\end{equation}
 with the \ $a_{2p}$\ obeying the relation
\begin{equation}
({\pi\over 2})^2 [(2p+2)(2p+1)a_{2p+2} - 4p^2 a_{2p}] - 12 e_3 a_{2p} + 
6 \sum_{0<r<p} a_{2r} a_{2p-2r} = 0 
\end{equation} 
where  \ $ e_3 = - \frac {(\pi)^2}{12} [\theta_2^4(0,\tau) + 
\theta_3^4(0,\tau)] $\ and \ $a_2 = \frac {\pi ^2}{4} \theta_2^4(0,\tau)  
\theta_3^4(0,\tau) .$\\
Moreover, this expansion is valid for $Im \tau > 0$ and for $\mid 
\sin\pi v \mid < 1  . $}

\bigskip
After eliminating the \ $a_{2p}$ \ we get relations between \ $c_{2p}$\ 
coefficients only. \\
The next result yields an analog version of the preceding 

\bigskip
	
	{\bf Proposition 2-4 }\quad {\it The theta function\ 
$\theta_4(v,\tau)$\  may be expressed under the form $$\theta_4(v,\tau) = 
\theta_4(0,\tau)\ \exp[ \sum_{p\geq 1}
	c_{2p}(\tau) (\sin\pi v)^{2p}]$$
	  where coefficients \ $c_{2p} $ \ satisfy the system of recurrence 
equations for \ $p \geq 1$ 
	$$ (A)
	 \quad \cases{
	4 ! \ {2p+4 \choose 4}\ c_{2p+4}  =  (2p+1)(2p+2) \big [(2p+2)(2p+3) + 
4 p^2 - c_0 \big ]c_{2p+2} & \cr
	  + (2p)^2 [ c_0 - (2p)^2 ] c_{2p} - 6 \big [ (2p+1)(2p+2) c_{2p+2} - 
2c_2 - \sum_{k=1}^p 2k c_{2k} \big ]^2  & \cr}$$
	where \quad  $ c_0 = - 4 [\theta_2^4(0,\tau) + \theta_3^4(0,\tau)], 
\quad  c_2 = \frac {1}{2 \pi ^2}  \frac {\theta ''_4 (0,\tau )}{\theta _4 
(0,\tau ) } $  and \\ $ c_4 =  \frac {1}{3} \theta_2^4(0,\tau)  
\theta_3^4(0,\tau) + \frac {1}{3} c_2.$\\
	Moreover, this expansion is valid for $Im \tau > 0$ and for $\mid 
\sin\pi v \mid < 1  . $}
	
	\bigskip

Indeed, starting from (2) a calculation gives 
$$\cases{
\frac {4}{\pi^2} [(2p+2)(2p+1)a_{2p+2} - 4p^2 a_{2p}] =  
(2p+1)(2p+2)(2p+3)(2p+4)c_{2p+4}  - & \cr 
\qquad (2p+1)(2p+2) [(2p+2)(2p+3) + 4 p^2]c_{2p+2}  + (2p)^4 c_{2p} & 
\cr}. $$
We also have
$$\sum_{0<r<p} a_{2r} a_{2p-2r} = (\frac {\pi }{2})^4 [ (2p+1)(2p+2) 
c_{2p+2} - c_2 - \sum_{k=1}^{k=p} c_{2k}]^2.$$
Thus, we obtain the following expressions connecting the $c_{2p}$ 
coefficients
$$ \quad \cases{
(2p+1)(2p+2)(2p+3)(2p+4)c_{2p+4}  - (2p+1)(2p+2) \big [(2p+2)(2p+3)  & 
\cr 
+ 4 p^2 \big ]c_{2p+2}  + (2p)^4 c_{2p} + c_0  \big [(2p+1)(2p+2) 
c_{2p+2} \ - & \cr 
(2p)^2 c_{2p} \big ] + 6 \big [ (2p+1)(2p+2) c_{2p+2} - 2c_2 - 
\sum_{k=1}^p 2k c_{2k} \big ]^2 = 0 & \cr}$$
which precisely are the recursion relations (A).
\subsection{Main result}
 
Now we are going to prove the following

\bigskip 
 {\bf Theorem 2-5} \quad {\it Let \ $q = e^{i \pi \tau }, \quad \mid 
q\mid < 1.$\ The coefficients \ $c_{2p}$\ defined by Propostion 2-4 may 
be expressed as
$$c_{2p}(\tau ) = - \frac {1}{p} \sum_{k\geq 0} \frac {1}{(\sin 
(k+\frac {1}{2})\pi \tau)^{2p}} = - \frac {1}{p} \sum_{k\geq 0} \bigg [\frac 
{(-4) q^{2k+1} }{(1 - q^{2k+1})^2}\bigg ]^p.$$
Then, the theta functions have the following expansions 
$$\theta_4(v,\tau) = \theta_4(0,\tau)\ \exp[- \sum_{p\geq 1} \sum_{ 
k\geq 0} \frac {1}{p} \bigg( \frac {\sin \pi v}{(\sin (k+\frac {1}{2})\pi 
\tau)}\bigg)^{2p}],$$
$$\theta_3(v,\tau) = \theta_4(0,\tau)\ \exp[- \sum_{p\geq 1} \sum_{ 
k\geq 0}\frac {1}{p} \bigg( \frac {\cos \pi v}{(\sin (k+\frac {1}{2})\pi 
\tau)}\bigg)^{2p}],$$
$$\theta_2(v,\tau) = \theta_4(0,\tau)\ \exp[ i\pi (v+\frac {1}{4}\tau)- 
\sum_{p\geq 1} \sum_{ k\geq 0}\frac {1}{p} \bigg( \frac {\cos \pi 
(v+\frac {1}{2}\tau)}{(\sin (k+\frac {1}{2})\pi \tau)}\bigg)^{2p}],$$
$$\theta_1(v,\tau) = \theta_4(0,\tau)\ \exp[ i\pi (v-\frac {1}{2}+\frac 
{1}{4}\tau) - \sum_{p\geq 1} \sum_{k\geq 0}\frac {1}{p} \bigg( \frac 
{\sin \pi (v+\frac {1}{2}\tau)}{(\sin (k+\frac {1}{2})\pi 
\tau)}\bigg)^{2p}].$$ 
Moreover, the above expressions for \ $\theta_4 $ \ and \ $ \theta_3$ \ 
are valid in the "strip" \ $\mid \frac {\sin \pi v}{\sin \frac 
{1}{2}\pi \tau} \mid < 1    ,$\\ those for  \ $\theta_2 $ \ and \ $ \theta_1$ \  
are valid in the "strip" \ $\mid \frac {\sin \pi (v+\frac {1}{2}\tau 
)}{\sin \frac {1}{2}\pi \tau} \mid < 1  .$}\\

\bigskip

{\bf Proof of Theorem 2-5}\quad  
At first, let us  determine  the coefficient \ $c_2(\tau ).$  

\bigskip

{\bf Lemma  2-6} \quad {\it The coefficient \ $c_2(\tau )$\ of the 
expansion of theta function given by Proposition 2-4 may be written under 
the forms }
$$c_2(\tau ) = - 4 \frac {{\displaystyle \sum_{n\geq 0} (-1)^n n^2 
q^{n^2}}} {{\displaystyle 1 + 2 \sum_{n\geq 1} (-1)^n q^{n^2}}}=  4 
\sum_{n\geq 1} \frac {q^{2n-1}}{(1 - q^{2n-1})^2} = 4 \sum_{n\geq 1} \frac {n 
q^{n}}{(1 - q^{2n})}$$
 
\bigskip

Indeed, from Propostion 2-4 one has $$2 \pi ^2 c_2 = \frac {\theta 
_4''(0)}{\theta _4(0)} .$$
Furthermore, by the Fourier series expansion $$\theta _4(v) = 1 + 2 
\sum_{n\geq 0} (-1)^n e^{i n^2 \pi \tau }\cos(2n \pi v),$$
 one deduces \ $\theta _4(0) = \sum_{n\geq 1} (-1)^n e^{i n^2 \pi \tau 
} $\ and 
$$\theta _4''(0) = - 8 \pi ^2 \sum_{n\geq 1} (-1)^n n^2 e^{i n^2 \pi 
\tau } .$$
On the other hand, since we have (see [W-W] or [B] vol 3)
$$\frac {\theta _4''(0)}{\theta _4(0)}  = 8 \pi ^2 \sum_{n\geq 1} \frac 
{q^{2n-1}}{(1 - q^{2n-1})^2 } $$ 
we then get the second expression of \ $c_2,$\ where \ $q = e^{i \pi 
\tau }$.

\bigskip

We shall compare different expansions for the theta functions. This 
will allow us to establish the expression of the $c_{2p}$ coefficients. \\
By Proposition 2-2
\begin{equation}
 \theta_4(v,\tau) = \theta_4(0,\tau)\ \exp[4 \sum_{n\geq1} 
{q^n\over{1-q^{2n}}} {(\sin{n\pi 
v})^2\over{n}}].
\end{equation}

 Comparing the latter with the expansion given by Proposition 2-4
$$\theta_4(v,\tau) = \theta_4(0,\tau)\ exp[ \sum_{p\geq 1}
c_{2p}(\tau) (\sin\pi v)^{2p}]$$
one derives the following

\bigskip

{\bf Lemma 2-7} \quad {\it The coefficients \ $c_{2p}$\ defined by 
Proposition 2-4 may be written as
$$c_{2p}(\tau ) =  (-1)^{p+1} 2^{2p} \bigg[ \frac {q^p}{p(1 - q^{2p})} 
+  2 \sum_{m > p} {m+p-1 \choose m-p-1}  \frac {q^m}{(m-p)(1 - 
q^{2m})}\bigg].$$
$$= (-1)^{p+1} \frac {2^{2p+1}}{(2p)!} \sum_{n \geq p} \frac 
{(n+p-1)!}{ (n-p)!}  \frac {q^n}{(1 - q^{2n})}.$$ 
 In particular, we again find }
$$c_2(\tau ) = 4 \sum_{n\geq 1} \frac {n q^{n}}{(1 - q^{2n})}.$$

\bigskip

{\bf Proof of Lemma 2-7} \quad It is enough to express $\cos (2n\alpha 
)$ in terms of $\sin \alpha$ as
$$ \cos (2n\alpha ) = \frac {(-1)^n}{2} (2 \sin \alpha )^{2n} + 
\sum_{k=0}^{n-1} (-1)^{n+k+1} \frac {2n}{k+1} {2n-k+1\choose k} (2\sin \alpha 
)^{2n-2k-2}$$ 
$$= n \sum_{0 \leq p \leq n} \frac {(-1)^p (n+p-1)!}{ (2p)! (n-p)!} 
(2\sin x)^{2p} = {\displaystyle \sum_{p=0}^{n}} b_{2p,n}(\sin \alpha 
)^{2p},$$
 where 
$$b_{2n,n} = (-1)^n 2^{2n-1}; \quad {\it and}\quad b_{2p,n} = (-1)^p 
2^{2p} \frac {n}{n-p} {n+p-1\choose n-p-1}\ {\it for } \ n \neq p.$$
It can also be rewritten
$$b_{2p,n} = n \frac {(-1)^p (n+p-1)!}{ (2p)! (n-p)!} 2^{2p}\quad {\it 
for}\quad 0 \leq p \leq n.$$
Furthermore, we may use the fact  $$ \log \frac {\theta _4 (v,\tau 
)}{\theta _4 (0,\tau )} = 2 \sum_{n\geq 1} \gamma _n - 2 \sum_{n\geq 1} 
\gamma _n \cos (2n v ) =  \sum_{p\geq 1} c_{2p}(\tau) (\sin\pi v)^{2p}, $$
where \ ${\displaystyle \gamma _n = \frac {q^n}{n (1-q^{2n})}}.$ \\
So, we obtain the relations
$$c_{2p} = - 2 \sum_{n\geq p} \gamma _n b_{2p,n} = - (-1)^p 
2^{2p}\bigg[ \gamma _p b_{2p,p} + \sum_{n\geq p} \frac {q^n}{(n-p)(1-q^{2n}} 
{n+p-1\choose n-p-1}\bigg].$$
Therefore 
$$c_{2p} = - 2 (-1)^p 2^{2p-1} \frac {q^p}{1 - q^{2p}} + (-1)^{p+1} 
2^{2p+1}  \sum_{m > p} {m+p-1 \choose m-p-1}  \frac {q^m}{(m-p)(1 - 
q^{2m})}.$$ 
i.e., 
$$c_{2p}(\tau) = (-1)^{p+1} \frac {2^{2p+1}}{(2p)!} \sum_{n \geq p} 
\frac {(n+p-1)!}{ (n-p)!}  \frac {q^n}{(1 - q^{2n})}.$$
In particular, for \ $p = 1,$  $$c_{2}(\tau) =  2^{2}  \frac {q}{1 - 
q^{2}} +  8 \sum_{m > 1} {m \choose m-2}  \frac {q^m}{(m-1)(1 - 
q^{2m})},$$ and 
$$c_2(\tau ) = 2^{2}  \frac {q}{1 - q^{2}} +  4 \sum_{m > 1}   \frac {m 
q^m}{(1 - q^{2m})} = 4 \sum_{n\geq 1} \frac {n q^{n}}{(1 - q^{2n})}.$$
These series can be expressed as
$$(-1)^{p+1}\frac {(2p)!}{2^{2p+1}}c_{2p}(\tau ) = \sum_{n \geq p} 
\frac {(n+p-1)!}{ (n-p)!}  \frac {q^n}{(1 - q^{2n})} = \sum_{n \geq p} 
\frac {q^n(n+p-1)!}{ (n-p)!} \sum_{k \geq 0} q^{2kn}$$
$$= \sum_{k \geq 0,n\geq p} \frac {(n+p-1)!}{ (n-p)!} q^{2kn+n} =  
\sum_{k \geq 0} q^{(2k+1)p}
 \sum_{n \geq p} \frac {(n+p-1)!}{ (n-p)!} q^{(2k+1)(n-p)} $$ 
$$= \sum_{k \geq 0} g_p (q^{2k+1}),$$
where 

$$g_p (z) = \sum_{n \geq p} \frac {(n+p-1)!}{ (n-p)!} z^n = z^p 
\sum_{m\geq 0} \frac {(m+2p-1)!}{ (m)!} q^m.$$
Now using the notation 
$$(2p)_n = 2p (2p + 1).....(2p + n - 1),\quad {\it then}\quad (m + 2p 
-1) ! = (2p)_m (2p - 1) !,$$
and one gets $$\frac {1}{(2p - 1) !} g_p (z) = z^p \sum_{m\geq 0} \frac 
{(2p)_m}{ (m)!} q^m.$$
However, we know that $$z^p \sum_{m\geq 0} \frac {(2p)_m}{ (m)!} q^m \ 
= z^p \ _2F_1 (2p,\alpha,\alpha,z) = \frac {z^p}{(1 - z)^{2p}}.$$
 $$g_p(\tau) = (2p-1)! \frac {z^p}{(1 - z)^{2p}}.$$
As a result, we may derive the expression 
$$c_{2p}(\tau ) = (-1)^{p+1} \frac {2^{2p+1}}{(2p)!} \sum_{k \geq 0} 
g_p (q^{2k+1}) = (-1)^{p+1} \frac {2^{2p+1}}{2p} \sum_{k \geq 0} \frac 
{q^{(2k+1)p}}{(1 - q^{2k+1})^{2p}}.$$
$$(-1)^{p+1}\frac {(2p)!}{2^{2p+1}}c_{2p}(\tau ) = (2p-1)!\sum_{k\geq 
0} \ _2F_1 (2p,\alpha,\alpha,q^{2k+1})$$
Finally,
$$c_{2p}(\tau ) = (-1)^{p+1} \frac {2^{2p+1}}{2p} \sum_{k \geq 0}\bigg( 
\frac {q^{2k+1}} {(1 - q^{2k+1})^{2}}\bigg)^p = - \frac {1}{p} 
\sum_{k\geq 0} \frac {1}{(\sin (k+\frac {1}{2})\pi \tau)^{2p}}.$$
The proof of Theorem 2-5 has been achieved.

\subsection {The link with elliptic functions}

Now, consider the zeta function of Jacobi. It is defined by 
$$Zn(z,k) =  \frac {1}{2K} \frac {d}{dz} \log \theta _4(v,\tau ),$$
where \ $v = \frac {z}{2K}$\ and \ $K = 2\int_0^{\frac {\pi }{2}}\frac 
{dx}{\sqrt {1 - k^2 \sin ^2x}}$\ is the complete elliptic integral of 
the first kind and the modulus is such that \ $0 < k < 1$.\\
We have

\bigskip
 
 {\bf Corollary 2-8} \quad {\it  The zeta function of Jacobi has the 
following form
$$Zn(z,k) = \frac {\pi}{2K} \sin (\pi 2v) \sum_{k\geq 0} \frac {1}{\sin 
^2(\pi v) - \sin ^2(k+\frac {1}{2}\pi \tau)}$$
where \ $v = \frac {z}{2K}$ satisfies \ $\mid  {\sin \pi v}\mid < 
\mid{(\sin (\frac {1}{2})\pi \tau)} \mid .$\\
In particular, the logarithmic derivatives of theta functions can be 
written under the forms}
$$\frac {\theta ' _4(v,\tau )}{\theta _4(v,\tau )} = 4 {\pi} \sin (\pi 
2v) \sum_{k\geq 0} \frac {q^{2k+1}}{1 - 2 q^{2k+1} \cos 2\pi v + 
q^{4k+2}}$$
$$\frac {\theta ' _3(v,\tau )}{\theta _3(v,\tau )} = - 4 {\pi} \sin 
(\pi 2v) \sum_{k\geq 0} \frac {q^{2k+1}}{1 + 2 q^{2k+1} \cos 2\pi v + 
q^{4k+2}}$$
$$\frac {\theta ' _2(v,\tau )}{\theta _2(v,\tau )} = - \tan (\pi v) - 4 
{\pi} \sin (\pi 2v) \sum_{k\geq 0} \frac {q^{2k+2}}{1 + 2 q^{2k+2} \cos 
2\pi v + q^{4k+4}}$$
$$\frac {\theta ' _1(v,\tau )}{\theta _1(v,\tau )} = \cot (\pi v) + 4 
{\pi} \sin (\pi 2v) \sum_{k\geq 0} \frac {q^{2k+2}}{1 - 2 q^{2k+2} \cos 
2\pi v + q^{4k+4}}.$$
{\it Moreover, the equations for \ $\theta _1$\ and \ $\theta _2$\ are 
valid in the strip  \ $\mid Im v \mid < Im \tau  ,$\  those for  \ 
$\theta _3$\ and \ $\theta _4$\ are valid in the strip} \ $\mid Im v \mid < 
\frac {1}{2} Im \tau  .$

\bigskip

Indeed, 
$$Zn(z,k) =  \frac {1}{2K} \frac {d}{dz} \log \theta _4(v,\tau ) = 
\frac {\pi}{2K} \sin (\pi 2v) \sum_{k\geq 0} \sum_{p\geq 1}  \bigg( \frac 
{\sin \pi v}{(\sin (k+\frac {1}{2})\pi \tau)}\bigg)^{2p}.$$

Suppose the variable $v$ satisfies
$$\mid \frac {\sin \pi v}{(\sin (\frac {1}{2})\pi \tau)} \mid < 1.$$

We then obtain
$$ \bigg (\frac {\sin \pi v}{(\sin (k+\frac {1}{2})\pi \tau)}\bigg )^2 
\sum_{p\geq 0} \bigg (\frac {\sin \pi v}{(\sin (k+\frac {1}{2})\pi 
\tau)}\bigg )^{2p} = \frac { \bigg (\frac {\sin \pi v}{(\sin (k+\frac 
{1}{2})\pi \tau)}\bigg )^2 }{1 -     
 \bigg (\frac {\sin \pi v}{(\sin (k+\frac {1}{2})\pi \tau)}\bigg )^2} = 
$$
 $$  \frac { ( \sin \pi v )^2}{{(\sin (k+\frac {1}{2})\pi \tau)}^2 -  
(\sin \pi v )^2}.$$
Therefore, the result follows. The domain of convergence for these 
series may be extended to the strip \ $\mid Im v \mid < \frac {1}{2} Im 
\tau $\ (see for example [W-W] page 489). \\
Notice that the zeta function of Jacobi also has a Fourier expansion 
$$Zn(z,k) = \frac {2\pi}{K} \sum_{n\geq 1}\frac {q^n}{1-q^{2n}} \sin 
\frac {n \pi z}{K} .$$

\bigskip

{\bf Remark  2-10} \quad  {\bf (i)}\quad Concerning the convergence of 
the trigonometric series, notice that those pertainig to  $\theta_3$ 
and $\theta_3$ could converge if $cosh (Im v) < \sinh \frac {Im 
\tau}{2}$. Those applicable to $\theta_1$ and $\theta_2$ could converge if $\mid 
\coth (Im v) \mid < \tanh \frac {Im \tau}{2}$.\\ The first inequality 
implies  $\mid \frac {\sin \pi v}{\sin \frac {1}{2}\pi \tau} \mid < 1$ 
and the second one implies  $\mid \frac {\sin \pi (v+\frac {1}{2}\tau 
)}{\sin \frac {1}{2}\pi \tau} \mid < 1.$

{\bf (ii)}\quad Moreover,  the above trigonometric expansions of theta 
functions which are very closed to Fourier series expansions seem to be 
new. We have not yet found  an analog of these series (in this present 
form) nor an allusion to them in the classical literature concerning 
the theta functions. The only thing which may catch our attention is the 
connection of the above expansions of theta functions with the infinite 
products. \\
Indeed, The latter expressions given by Corollary 2-8 can also be 
deduced by logarihmic differentiation of the wellknown identities which 
yields the theta functions as infinite product. Therefore, the series 
expansions stated by Theorem 2-5 can also be deduced.\\
More precisely, from
$$\theta _1 (v,\tau) = 2 q^{\frac {1}{4}} \sin (\pi z) \Pi _{k\geq 0} 
(1 - q^{2k+2}) (1 - 2 q^{2k+2} \cos (2\pi v) + q^{4k+4})$$ 
$$\theta _2 (v,\tau) = 2 q^{\frac {1}{4}} \cos (\pi z) \Pi _{k\geq 0} 
(1 - q^{2k+2}) (1 + 2 q^{2k+2} \cos (2\pi v) + q^{4k+4})$$
$$\theta _3 (v,\tau) = \Pi _{k\geq 0} (1 - q^{2k+2}) (1 + 2 q^{2k+1} 
\cos (2\pi v) + q^{4k+2})$$
$$\theta _4 (v,\tau) = \Pi _{k\geq 0} (1 - q^{2k+2}) (1 - 2 q^{2k+1} 
\cos (2\pi v) + q^{4k+2})$$ 
we get the above expression of \ ${\displaystyle \frac {\theta ' 
_1(v,\tau )}{\theta _1(v,\tau )}}$, ${\displaystyle \frac {\theta ' _2(v,\tau 
)}{\theta _2(v,\tau )}}$, ${\displaystyle \frac {\theta ' _3(v,\tau 
)}{\theta _3(v,\tau )}}$, ${\displaystyle \frac {\theta ' _4(v,\tau 
)}{\theta _4(v,\tau )}}$,\  respectively (see [W-W] page 489).\\ However, we 
may notice that the previous infinite product representing the theta 
functions are valid in the entire  $v$-plan. 
  
\bigskip

{\bf Corollary 2-11}\quad
{\it Under the same hypotheses, the following expressions for ratios of 
theta functions hold}
$$\frac {\theta_1(v,\tau)}{\theta_2(v,\tau)} = exp [\sum_{p\geq 1} 
c_{2p}(\tau) [ \sin^{2p}\pi (v+{1\over 2}\tau) - \cos^{2p}\pi (v+{1\over 
2}\tau)] ]$$
$$\frac {\theta_3(v,\tau)}{\theta_4(v,\tau)} = exp [\sum_{p\geq 1} 
c_{2p}(\tau) [ \cos^{2p}\pi v - \sin^{2p}\pi v] ] .$$

\bigskip

{\bf Corollary 2-12}\quad
{\it Under the same hypotheses, the following expression for the 
product of the theta functions holds}
$$\frac {\theta_2(v,\tau) \theta_3(v,\tau) 
\theta_4(v,\tau)}{\theta_4^3(0,\tau)} = e^{i\pi (v+\frac {1}{4}\tau)} exp [ \sum_{p\geq 1} 
c_{2p}(\tau) [ \cos^{2p}\pi v + \sin^{2p}\pi v+ \cos^{2p}\pi (v+\frac 
{\tau}{2})] ] .$$
In particular, we get 
$$\theta_1'(0,\tau) = \pi \theta_2(0,\tau) \theta_3(0,\tau) 
\theta_4(0,\tau) = -\pi \theta_4^3(0,\tau) q^{1\over 4} exp [\sum_{p\geq 1} 
c_{2p}(\tau) [ 1 + \cos^{2p}\pi {\tau \over 2}] ].$$
The eta function of Dedekind may be expressed as
$$\eta (\tau) = 2^{-\frac {1}{3}} e^{i\pi \tau /12} \theta_4(0,\tau )\ 
exp \sum_{p\geq 1} \frac {c_{2p}(\tau)}{3} [1 + \cos^{2p}(\pi \frac 
{\tau}{2})]$$

The proofs of the above results use the same techniques of [C].
\newpage
\section{Modular transformations}

\subsection{Some properties of the coefficients}
The above coefficients \ $c_{2p}(\tau) $\ are such that 
$$c_{2p}(\tau+2) = c_{2p}(\tau).$$
Furthermore, by the periodicity properties of the theta functions we 
may deduce some additional relations.
More precisely,  if the argument $v$  is increased by $\tau$, they are 
unaffected except for the multiplication by a simple factor. We get 
$$\theta _1 (v+\tau,\tau) = - (q e^{2i\pi v})^{-1} \theta _1 (v,\tau) 
$$
$$\theta _2 (v+\tau,\tau) =  (q e^{2i\pi v})^{-1} \theta _2 (v,\tau) $$
$$\theta _3 (v+\tau,\tau) =  (q e^{2i\pi v})^{-1} \theta _3 (v,\tau) $$
$$\theta _4 (v+\tau,\tau) = - (q e^{2i\pi v})^{-1} \theta _4 (v,\tau) 
$$
We deduce the following
$$\sum_{p\geq 1} c_{2p}(\tau) [(\cos \pi v)^{2p} - (\sin \pi v)^{2p}] = 
i \pi + \sum_{p\geq 1} c_{2p}(\tau) [(\cos \pi (v+\tau))^{2p} - (\sin 
\pi (v+\tau))^{2p}]. $$
$$\sum_{p\geq 1} c_{2p}(\tau) [(\cos \pi v)^{2p} - (\cos \pi 
(v+\tau))^{2p}] = i\pi (2v + \tau)$$
 \subsection{Actions of the modular group}

Moreover, we have seen before (see [C])  that the double family of  
coefficients \ $a_{2p}$ \ and \ $c_{2p}$,
have analogous properties under the action of the modular group
 \ $\Gamma(1)$. For example, when we examine the theta relation \ 
$\theta_4(v+{1\over 2},\tau+1) =
\theta_4(v,\tau)$,\ we find

\begin{equation}
c_{2p} (\tau+1) = (-1)^p  \sum_{k\geq p}  {k \choose p} c_{2k}(\tau) 
\end{equation}
which is similar with the relation satisfied by the coefficients  \ 
$a_{2p}$\ (see [C]).\\
Equality (5), combined with System (A), permits to obtain other 
relations between the coefficients.\\
As consequences of Theorem 2-5, we have
\bigskip

{\bf Corollary 3-1} \quad {\it Under the actions of the modular group \ 
$\tau \rightarrow \tau + 1,$\ and \ $\tau \rightarrow \frac {-1}{\tau } 
,$\ the coefficients become
$$c_{2p}(\tau +1) =  - \frac {1}{p} \sum_{k\geq 0} \bigg [\frac {4 
q^{2k+1} }{(1 + q^{2k+1})^2}\bigg ]^p.$$
$$\frac {2}{\tau ^2}c_{2p}(\frac {-1}{\tau}) = -1 - (-1)^p\frac 
{2^{p+1}}{2p} \sum_{k\geq 1}\bigg [\frac {q^{2k}}{1+q^{2k}}\bigg ]^p.$$
Moreover, the following relation holds}
$$(-1)  \sum_{k\geq 1}  k c_{2k}(\tau) =  \sum_{n\geq 1} \frac {n 
(-1)^n q^{n}}{(1 - q^{2n})}.$$
$$\sum_{p\geq 1} c_{2p}(\tau) [(\cos \pi \frac {\tau}{2})^{2p} = 
\sum_{p\geq 1} c_{2p}(\tau+1) [1 + (\sin \pi \frac {\tau}{2})^{2p}].$$
$$c_{2}(\tau +1) =  - 4 \sum_{k\geq 0} \frac { q^{2k+1} }{(1 + 
q^{2k+1})^2} = \frac {1}{2\pi ^2}\frac {\theta_3 ''(0)}{\theta_3(0)} .$$
$$\frac {2}{\tau ^2}c_{2}(\frac {-1}{\tau}) = -1 - 8 \sum_{k\geq 
1}\frac {q^{2n}}{1+q^{2n}} = \frac {1}{2\pi ^2}\frac {\theta_2 
''(0)}{\theta_2(0)}.$$

\bigskip
\subsection{Landen and Gauss transformations}

Furthermore, using the Landen transformation we may obtain relations 
connecting \ $\theta (2v,2\tau )$\ and \ $\theta (v,\tau ).$\
 More precisely, we get (see [B] or [W-W])
\begin{equation}
\theta _4(2v,2\tau ) = \frac {\theta _3(v,\tau ) \theta _4(v,\tau 
)}{\theta _4(0,2\tau )} = \theta_4(0,2\tau)\ \exp[ \sum_{p\geq 1}
c_{2p}(2\tau) (\sin 2\pi v)^{2p}].
\end{equation} 
From Proposition 2-4 and Theorem 2-5 , one proves the following

\bigskip

{\bf Proposition 3-2} \quad  {\it The coefficients \ $c_{2p}(2\tau )$ \ 
satisfy the relation}
$$c_{2p}(2\tau +1) =  \sum_{k\geq 2p} 2^{-k} {k \choose 2p} c_{2k}(\tau 
) $$

\bigskip 
Indeed, we have seen that\ 
$$\theta _3(v,\tau ) = \theta_4(0,\tau)\ \exp[ \sum_{p\geq 1}
c_{2p}(\tau) (\cos \pi v)^{2p}].$$ 
Moreover, since \ $[\theta _4(0,2\tau )]^2 = \theta _3(0,\tau ) \theta 
_4(0,\tau )$ \ and \ $ \log[\theta _3(0,\tau )] = \sum_p c_{2p}(\tau 
)$, 
then 
(6) implies\\
 
\bigskip
{\bf Lemma 3-3}\quad {\it When $v$ belongs to the strip \ $\mid \frac 
{\sin \pi v}{\sin \frac {1}{2}\pi \tau} \mid < 1    $\  the following 
relations hold
$$
 \sum_{p\geq 1} c_{2p}(\tau) \bigg[ (\sin \pi v)^{2p} + (\cos \pi 
v)^{2p} - 1\bigg] = \sum_{p\geq 1}
c_{2p}(2\tau) (\sin 2\pi v)^{2p}$$ 
 $$ = \sum_{p\geq 1}
c_{2p}(2\tau +1) (\cos 2\pi v)^{2p}.$$ 
$$
 \sum_{p\geq 1} c_{2p}(\tau) \bigg[ (\sin \pi (v+\frac {1}{4}))^{2p} + 
(\cos \pi (v+\frac {1}{4}))^{2p} - 1\bigg] = \sum_{p\geq 1}
c_{2p}(2\tau) (\cos 2\pi v)^{2p}$$ 
 $$ = \sum_{p\geq 1}
c_{2p}(2\tau +1) (\sin 2\pi v)^{2p}.$$
In particular, we deduce the identities}
$$ \theta _3(2v,2\tau ) = \frac {\theta _3(v+\frac {1}{4},\tau ) \theta 
_3(v+\frac {1}{4},\tau )}{\theta _4(0,2\tau )} = \frac {\theta 
_4(v+\frac {1}{4},\tau ) \theta _4(v+\frac {1}{4},\tau )}{\theta _4(0,2\tau 
)}$$

\bigskip

Indeed, since \ $ (\cos v)^{2n} = ( \frac {1 + \cos 2v}{2})^n = 2^{-n} 
{\displaystyle \sum_{k=0}^{k=n}} {n\choose k} (\cos 2v)^k$ \ then 
$$\sum_p 2^{-p} c_{2p}(\tau)\sum_{k=0}^{k=p}[ (-1)^k {p\choose k}+{p 
\choose k}] (\cos 2v)^k-\sum c_{2p}(\tau)=$$
 $$ = \sum c_{2p}(2\tau +1) (\cos 2v)^{2p},$$ thus 
$$ \sum_p 2^{-p} c_{2p}(\tau) \sum_{k=0}^{2k=p}  {p\choose 2k} (\cos 
2v)^{2k} -  \sum c_{2p}(\tau) = \sum c_{2p}(2\tau +1) (\cos 2v)^{2p}.$$
As a result 
$$ \sum_k \bigg[ \sum_{p\geq 2k} 2^{-p} {p\choose 2k}c_{2p} (\tau) 
\bigg]  (\cos 2v)^{2k} -  \sum c_{2p}(\tau) = \sum c_{2p}(2\tau +1) (\cos 
2v)^{2p}.$$
Interchanging \ $p$ \ and \ $k,$\ one gets
$$ \sum_p \bigg[ \sum_{k\geq 2p} 2^{-k} {p\choose 2k}c_{2k} (\tau) 
\bigg]  (\cos 2v)^{2p} -  \sum c_{2p}(\tau) = \sum c_{2p}(2\tau +1) (\cos 
2v)^{2p}.$$
This proves Lemma 3-3 and Proposition 3-2.

\bigskip

{\bf Corollary 3-4}\quad {\it The coefficients \ $c_{2p}(2\tau )$ \ 
also satisfy the equation}
$$c_{2p}(2\tau ) = (-1)^p  \sum_{m\geq 2p} [ \sum_{k=p}^{2k=m} 
{k\choose p}{m\choose 2k}] 2^{-m}  c_{2m}(\tau ) $$

\bigskip

Indeed, (5) implies \ $c_{2p} (\tau+1) = (-1)^p  \sum_{k\geq p}  {k 
\choose p} c_{2k}(\tau) $\ and \\ Proposition 3-2 gives \ $c_{2p}(2\tau 
+1) =  \sum_{k\geq 2p} 2^{-k} {k \choose 2p} c_{2k}(\tau ).$\\

Thus, combining these two equalities, one obtains
$$c_{2p}(2\tau ) = (-1)^p  \sum_{k\geq p}{k\choose p} [ \sum_{2k\geq 
m}{m\choose 2k}] 2^{-m}  c_{2m}(\tau ) $$
$$  = (-1)^p  \sum_{m\geq 2p} [ \sum_{k=p}^{2k=m} {k\choose p}{m\choose 
2k}] 2^{-m}  c_{2m}(\tau ) $$

\bigskip

{\bf Corollary 3-5}\quad {\it The $\eta $ function of Dedekind 
satisfies the following }
$$\frac {\eta (2\tau)}{\eta (\tau)} = 2^{-1/3} q^{\frac {1}{12}} exp 
\sum_{p\geq 1} \frac {c_{2p}(\tau)}{3} [\cos ^{2p}\frac {\pi \tau}{2} - 
1].$$

\bigskip
Indeed, this follows from the known identities 
$$\frac {\eta ^3(2\tau)}{\eta ^3(\tau)} = \frac {\theta 
'_1(0,2\tau)}{\theta '_1(0,\tau)} = \frac {1}{2} \frac {\theta 
_2(0,\tau)}{\sqrt{}\theta _3(0,\tau)\theta _4(0,\tau)}.$$

\newpage

\section{Some arithmetical properties} 

\subsection{Transformations of higher order}
One proved (Proposition (3-4) of [C1]) for \ $m \geq 1$ \ a modular 
equation satisfied by the coefficients
$$\frac {1}{2} \frac {q^m}{1 - q^{2m}} = m (-1)^{m+1} \sum_{p\geq m} 
2^{-2p} {2p\choose p-m} c_{2p}(\tau ) =  \sum_{p\geq 1} 2^{-2p} 
{2p\choose p-1} c_{2p}(m\tau ) .$$
In particular, for \ $m = 2$ \ one again finds  
$$- \sum_{p\geq 2} 2^{-2p} {2p\choose p-2} c_{2p}(\tau ) = \sum_{p\geq 
1} 2^{-2p} {2p\choose p-1} c_{2p}(2\tau )$$

In fact, for any positive integer \ $k$\ we get the identity
$$\frac {1}{2} \frac {q^m}{1 - q^{2m}} =\sum_{p\geq m}^{k'\geq 0} 
{2p\choose p-1} \frac {(-1)^{p+1}}{p} \bigg( \frac {q^{2k'm+m}} {(1 - 
q^{2k'm+m})^{2}}\bigg)^p $$ $$= k (-1)^{k+1} \sum_{p\geq k} 2^{-2p} 
{2p\choose p-k} c_{2p}(\frac {m}{k}\tau ).$$
Thus, we have

\bigskip

{\bf Proposition 4-1}\quad {\it Let \ $n \ {\it and}\ k \geq 1$\ be  
integers.\\
The following identities hold
$$k (-1)^{k+1} \sum_{p\geq k} 2^{-2p} {2p\choose p-k} c_{2p}( n\tau ) = 
n (-1)^{n+1} \sum_{p\geq n} 2^{-2p} {2p\choose p-n} c_{2p}( k\tau ) =$$
$$\sum_{p\geq 1} 2^{-2p} {2p\choose p-1} c_{2p}(nk\tau ) = nk 
(-1)^{nk+1} \sum_{p\geq nk} 2^{-2p} {2p\choose p-nk} c_{2p}(\tau )$$

\bigskip

\subsection{Other expressions using Lambert series}

{\bf Propostion 4-2}\quad {\it The coefficients \ $c_{2p}(\tau )$ \ 
given by Theorem 2-5 may also be written under the form 
$$c_{2p}(\tau ) = (-1)^{p+1} \frac {2^{2p+1}}{2p} \sum_{k \geq 0}\bigg( 
\frac {q^{2k+1}} {(1 - q^{2k+1})^{2}}\bigg)^p = (-1)^{p+1} \frac 
{2^{2p}}{ 2p\ !} \sum_{n\geq p} (A_n + B_n) q^{n}$$ where } }
$$ A_n + B_n = \sum_{d \mid n,d\geq p} (1 + (-1)^{\frac {n}{d} - 1}) 
\frac {(d + p - 1) !}{(d-p)!}$$ 

\bigskip

{\bf Proof}\quad It is known that for \ $\mid z\mid < 1$ \ Lambert 
series can be written in different manners. In particular, 
$$\sum_{n\geq 1} a_n \frac {z^n}{1 - z^n} = \sum_{n\geq 1} A_n z^n,$$
where  $$A_n = \sum_{d\mid n} a_d.$$ Moreover, using the Mobius 
function \ $\mu (n)$,\ we can express \ $a_n$\ in terms of \ $A_n$.\ Namely,
$$a_n = \sum_{d\mid n} \mu (\frac {n}{d})  A_n.$$
Furthermore, we may prove that 
$$\sum_{n\geq 1} \frac {z^n}{1 + z^n} = \sum_{n\geq 1} B_n z^n,$$
where  $$B_n = \sum_{t\mid n} (-1)^{t-1} a_{\frac {n}{t}} = \sum_{d\mid 
n} (-1)^{\frac {n}{d}}a_d.$$
Let for any \ $p \geq 1$ 
$$c_{2p}(\tau ) = (-1)^{p+1} \frac {2^{2p+1}}{ (2p) !} \sum_{n\geq p} 
\frac {(n + p - 1) !}{(n-p)!} \frac {q^n}{1 - q^{2n}}$$
Then, we may deduce
$$ 2 (-1)^{p+1} \frac {2^{2p}}{ 2p\ !}\sum_{n\geq p} \frac {(n + p - 1) 
!}{(n-p)!} \frac {z^n}{1 - z^{2n}} = \sum_{n\geq p} (A_n + B_n) z^n,$$
where $$A_n + B_n = \sum_{d \mid n,d\geq p} (1 + (-1)^{\frac {n}{d} - 
1}) \frac {(d + p - 1) !}{(d-p)!},$$
and \ $A_n = B_n = 0$\ for \ $n < p$. 

\bigskip
{\bf Remark 4-3} \quad  For \ $p = 1$, \ we get 
$$c_2 (\tau ) = 4 \sum_{n\geq 1} \frac {n q^{n}}{(1 - q^{2n})}$$
but, \ $$ \frac {n q^{n}}{(1 - q^{2n})} = n (q^n + q^{3n} + q^{5n} 
+....$$
 The coefficent of \ $q^m$\ is \ $n,$ \ if \ $n$\ divides \ $m$ \ in 
such a way that \ $\frac {m}{n} $\ is odd.\\ 
Suppose \ $ m = 2^k p$, \ $p$ being odd . It implies that \ $n = 2^k 
d,$\ where $d$ is a divisor of $p$.\\ Then, the coefficient of \ $q^m$\ 
in the expansion of \ $\frac {n q^{n}}{(1 - q^{2n})}$\ is obviously \ 
$2^k \sum _{d\mid p} d.$

\newpage 
\section{An expansion of \ $\wp(z)$} 

In this section we propose another interesting application of Theorem 
2-5, one concerning a new type of expansion for elliptic functions.
 More precisely, we state for the elliptic Weierstrass fonction 
$\wp(z)$ an analogous expansion as we did for the $\theta$ function.\\ 
   
Recall that $\wp(z)$  which has primitive periods \ $2$\ and \ $2\tau$ 
\ relative to \ $g_2$ \ and \ $g_3,$ may be written as
$$\wp (z) = \wp (z;1,\tau) = \frac {1}{z^2} + \sum_{m,n} \bigg [{\frac 
{1}{(z- 2m - 2n\tau)^2}} - {\frac {1}{(2m+2n\tau)^2}}\bigg ]$$
The Weierstrass fonction $\wp(z)$ is related to the theta functions 
$\theta_i (v)$ where $v = {z\over{2\omega}}$ :
$$\wp(z) = ({1\over{2\omega}})^2 [-4\eta \omega - {d^2 log 
\theta_1(v)\over{dv^2}}]$$ with $\eta = - \frac {1}{12\omega } \frac {\theta 
_1^{'''} (0)}{\theta _1^{'} (0)}$. In the same way we have 
$$
\wp(z+\omega ' ) = ({1\over{2\omega}})^2 [-4\eta \omega - {d^2 log 
\theta_4(v)\over{dv^2}}].
$$

We proved  (see [C] ) that \ $\wp(z;1,\tau)$\  has the following 
expansion
\begin{equation}
\wp(z+\tau) = e_3 - \sum_{p\geq 1} a_{2p}(e_3)  (\sin{\pi z\over 
{2}})^{2p}.
\end{equation}

where \ $a_2(e_3) = ({1\over{\pi}})^2 (g_2 - 12 {e_3}^2)$\ , \ $e_3 = 
\wp(\tau)$\ and the other coefficients satisfy the recursion relations
$$
 ({\pi\over{2}})^2 [(2p+2)(2p+1)a_{2p+2} - 4p^2 a_{2p}] - 12 e_3 a_{2p} 
+ 6 \sum_{0<r<p} a_{2r} a_{2p-2r} = 0.\quad (3)
$$
Using relations (2) and (3) and Theorem 2-5 we shall give explicit 
expressions for the coefficients $a_{2p}(e_3) = a_{2p}(\tau)$

\bigskip
  
{\bf Theorem 5-1} \quad {\it The Weierstrass elliptic function \ $\wp 
(z) = \wp (z;1,\tau)$\ with primitive periods $1$ and $\tau$ may be 
expressed under the form
$$\wp(z+\tau) = e_3 - \sum_{p\geq 1} a_{2p}(\tau)  (\sin{\pi z\over 
{2}})^{2p}$$
where the coefficients are
$$ a_{2p}(\tau ) = - 2(2p+1) \sum_{k\geq 0} \frac {1}{(\sin (k+\frac 
{1}{2})\pi \tau)^{2p+2}} + 4p \sum_{k\geq 0} \frac {1}{(\sin (k+\frac 
{1}{2})\pi \tau)^{2p}}.$$   
Moreover, the above expression of  $\wp (z)$ is valid in the strip  
$\mid \frac {\sin \pi v}{\sin \frac {1}{2}\pi \tau} \mid < 1 $.}\\

   \newpage
   Indeed, the expression for $ a_{2p}(\tau )$ is deduced from 
  $$(2p+2)(2p+1) c_{2p+2}(\tau) - 4p^2 c_{2p}(\tau) = \frac {4 
\omega^2}{ \pi ^2} a_{2p}(\omega,\omega') = \frac {3}{ \pi ^2} a_{2p}(1,\tau 
)\qquad (2)$$
  and  Theorem 4-1 since 
  $$c_{2p}(\tau ) = - \frac {1}{p} \sum_{k\geq 0} \frac {1}{(\sin 
(k+\frac {1}{2})\pi \tau)^{2p}} = - \frac {1}{p} \sum_{k\geq 0} \bigg [\frac 
{(-4) q^{2k+1} }{(1 - q^{2k+1})^2}\bigg ]^p.$$
  Thus, \ $\wp (z)$\ may be expressed as
  $$\wp(z+\tau) = e_3 - \sum_{p\geq 1}[  \sum_{k\geq 0} \frac {- 
2(2p+1)(\sin{\pi z\over {2}})^{2p}}{(\sin (k+\frac {1}{2})\pi \tau)^{2p+2}} +  
\sum_{k\geq 0} \frac {4p(\sin{\pi z\over {2}})^{2p}}{(\sin (k+\frac 
{1}{2})\pi \tau)^{2p}} ].$$
  Moreover, according to the addition theorem for \ $\wp(z),$\ this 
function may be written as
 $$\wp(z) = e_3 - \frac {g_2 - 12 e_3^2}{ {\displaystyle 2\sum_{p\geq 
1}  \sum_{k\geq 0}[ \frac {- 2(2p+1)(\sin{\pi z\over {2}})^{2p}}{(\sin 
(k+\frac {1}{2})\pi \tau)^{2p+2}} +  \frac {4p(\sin{\pi z\over 
{2}})^{2p}}{(\sin (k+\frac {1}{2})\pi \tau)^{2p}} ]}}.$$
  
\section{The Rogers-Ramanujan continued fractions}
For $\mid q\mid < 1$ the Rogers functions are defined by 
$$G(q) = \sum_{n\geq 0}\frac {q^{2n^2}}{(q^2,q^2)_n} ; \quad H(q) = 
\sum_{n\geq 0}\frac {q^{2n(n+1)}}{(q^2,q^2)_n}$$
where \ $q = e^{i\pi \tau}; \quad (a,q)_n = \prod_{k=0}^{n-1} (1 - a 
q^k); \quad (a,q)_\infty =  lim_{n\rightarrow \infty} (a,q)_n.$ \\ 
It is known that 
$$G(q) = \frac {1}{(q^2,q^{10})_\infty (q^8,q^{10})_\infty }; \quad 
H(q) = \frac {1}{(q^4,q^{10})_\infty (q^6,q^{10})_\infty }.$$
The Rogers-Ramanujan continued fractions is defined by
 $$R(q) = \frac {q^{2/5}}{1+\frac {q^2}{1+\frac {q^4}{1+\frac 
{q^6}{1+...}}}}.$$
It can be expressed in terms of $q$-series by 
$$ R(q) = q^{2/5} \frac {H(q)}{G(q)} =  q^{2/5} \prod_{k\geq 1} \frac 
{(1- q^{10k-2})(1- q^{10k-8})}{(1- q^{10k-4})(1- q^{10k-6})}$$
We intend to give a simple proof of the the above Rogers identity.

\bigskip 

{\bf Theorem 6-1} \quad {\it The Rogers-Ramanujan continued fractions \ 
$R(q)$ \ for \\ $q = e^{i\pi \tau}, \ \mid q \mid < 1$ \ may be 
expressed as 
$$ R(q) = exp [i\pi \frac {2\tau}{5} - \sum_{p\geq 1} c_{2p}(5\tau) 
[(\sin \frac {\pi \tau }{2})^{2p} - (\sin \frac {3\pi \tau }{2})^{2p}] ]$$
$$ = exp [i\pi \frac {2\tau}{5} + \sum_{p\geq 1}  \sum_{k\geq 0} + 
\frac {1}{p} \frac {[(\sin \frac {\pi \tau }{2})^{2p} - (\sin \frac {3\pi 
\tau }{2})^{2p}]}{(\sin (k+\frac {1}{2})5\pi \tau)^{2p}}].$$ 
Equivalently, in terms of $q$-series}
$$R(q) = q^{2/5} e^{ \sum_{p\geq 1,k\geq 0}{\displaystyle  
\frac{q^{2p(k+1)}}{p}\bigg(\frac {(q^2-q)^{2p}- 
(q^3-1)^{2p}}{(q^{5(2k+1)}-1)^{2p}}\bigg)}}.$$

\bigskip 
{\bf Proof }\quad Consider the Ramanujan theta function 
$$f(a,b) = \sum_{k\in Z} a^{\frac{k(k+1)}{2}} b^{\frac{k(k-1)}{2}}$$
$$f(a,b) = (-a,ab)_\infty (-b,ab)\infty (ab,ab)_\infty $$
with $\mid ab\mid < 1.$\ We denote here  \ $(\alpha,\beta )_\infty = 
\prod _{i\geq 1}(1 - \alpha \beta^i).$  
This function is related to the classical one
$$ \theta_4(v,q) = f(-qe^{2i\pi v},-qe^{-2i\pi v}).$$
Thus, $$R(q) = q^{1/5} \frac {f(-q^2,-q^8)}{f(-q^4,-q^6)}.$$
One gets the relation
$$R(q) = q^{2/5}\ \frac {\theta_4(\frac {3}{2}\tau, 
5\tau)}{\theta_4(\frac {1}{2}\tau, 5\tau)}, $$
which implies the wanted expression of $R(q)$ since by Theorem 2-5
$$\theta_4(v,5\tau) = \theta_4(0,5\tau)\ \exp[- \sum_{p\geq 1} \sum_{ 
k\geq 0} \frac {1}{p} \bigg( \frac {\sin \pi v}{(\sin (k+\frac 
{1}{2})5\pi \tau)}\bigg)^{2p}].$$

\bigskip 
{\bf Corollary 6-2} \quad {\it  $R(q)$ \  may also be expressed as    
   $$R(q) = q^{2/5} \exp \bigg( \sum_{p\geq 1,k\geq 0} \frac {1}{p} 
\frac {[(\cos \frac {\pi \tau }{2})^{2p} - (\cos \frac {3\pi \tau 
}{2})^{2p}]}{(\cos (k+\frac {1}{2})5\pi \tau)^{2p}}\bigg).$$ Equivalently,}
   $$R(q) =  q^{2/5} e^{ \sum_{p\geq 1,k\geq 0}{\displaystyle 
\frac{q^{2p(k+1)}}{p}\bigg(\frac {(q^2+q)^{2p}- 
(q^3+1)^{2p}}{(q^{5(2k+1)}+1)^{2p}}\bigg)}}.$$

\bigskip 
Indeed, this is a direct consequence of the theta identity
$$\theta_4(v,\tau +1) = \theta_3(v,\tau ) = \theta_4(v+\frac 
{1}{2},\tau).$$  

\bigskip 
{\bf Corollary 6-3} \quad {\it Under the action of the modular group 
one has}
$$R(\tau +1) = e^{2i\pi/5} R(\tau).$$
$$R(-1/\tau ) = e^{i\pi /5 + \sum_{p\geq 1,k\geq 0}{\displaystyle  
\frac {1}{p} \frac {[(\sin \frac {\pi (\tau -4 }{10})^{2p} - (\sin \frac 
{\pi (\tau -2)}{10})^{2p}]}{(\sin (k+\frac {1}{2})\pi \tau/5)^{2p}}}}.$$ 

\bigskip 

Indeed, the first identity is obtained from Corollary 6-2. Whereas, we 
may deduce the second one from
$$R(-1/\tau ) = 
e^{i\pi /5}\ \frac {\theta_3(\frac {3+\tau}{10}, 
\tau/5)}{\theta_3(\frac {1+\tau}{10}\tau, \tau/5)} =
 e^{i\pi /5}\ \frac {\theta_4(\frac {\tau -2}{10}, 
\tau/5)}{\theta_4(\frac {\tau -4}{10}\tau, \tau/5)} .$$ 
 
 \bigskip 
{\bf Corollary 6-4} \quad {\it Let $q = e^{i\pi \tau}, \ \mid q \mid < 
1.$ From the expression
$$R(q) = q^{2/5}\exp \bigg( \sum_{p\geq 1,k\geq 0} \frac {1}{p} \frac 
{[(\sin \frac {\pi \tau }{2})^{2p} - (\sin \frac {3\pi \tau 
}{2})^{2p}]}{(\sin (k+\frac {1}{2})5\pi \tau)^{2p}}\bigg)$$ 
 we may deduce the Rogers formula}
$$R(q) = q^{\frac {2}{5}} \prod_{k\geq 1} \frac {(1- q^{10k-2})(1- 
q^{10k-8})}{(1- q^{10k-4})(1- q^{10k-6})}.$$
 \bigskip

{\bf Proof} \quad Starting from the expansions
$$\sum_{p\geq 1} \frac {1}{p} \frac {(\sin \frac {\pi \tau 
}{2})^{2p}}{(\sin (k+\frac {1}{2})5\pi \tau)^{2p}} = - log (1 - \frac {(\sin \frac 
{\pi \tau }{2})^{2} }{(\sin (k+\frac {1}{2})5\pi \tau)^{2}}) $$
$$\sum_{p\geq 1} \frac {1}{p} \frac {(\sin \frac {3\pi \tau 
}{2})^{2p}}{(\sin (k+\frac {1}{2})5\pi \tau)^{2p}} = - log (1 - \frac {(\sin \frac 
{3\pi \tau }{2})^{2} }{(\sin (k+\frac {1}{2})5\pi \tau)^{2}}). $$ 
We then derive the following expressions
$$R(q) = q^{2/5} \prod_{k\geq 0} \frac {(\sin (k+\frac {1}{2})5\pi 
\tau)^{2} - (\sin \frac {3\pi \tau }{2})^{2}}{(\sin (k+\frac {1}{2})5\pi 
\tau)^{2} - (\sin \frac {\pi \tau }{2})^{2}}$$
$$R(q) =  q^{2/5} \prod_{k\geq 0} \frac {\cos (2k+1)5\pi \tau - \cos 
3\pi \tau}{\cos (2k+1)5\pi \tau - \cos \pi \tau}$$
$$R(q) = q^{2/5} \prod_{k\geq 0} \frac {[\sin (5k+4)\pi \tau ] [\sin 
(5k+1)\pi \tau ]}{[\sin (5k+3)\pi \tau ] [\sin (5k+2)\pi \tau ]}. $$
Thus,
$$R(q) = q^{2/5} \prod_{k\geq 0} \frac {(q^{10k+8} - 1)(q^{10k+2} - 
1)q^{5k+3} q^{5k+2}}{q^{5k+4} q^{5k+1}(q^{10k+6} - 1)(q^{10k+4} - 1)}.$$

{\bf Remark 6-5}\quad 
There are many other expressions of \ $R(q)$ \ in terms of theta 
functions. One of them given by Z.-G. Liu [L] may be interesting
$$R(q) = q^{-6/10} \frac 
{\theta_1(\tau,5\tau)}{\theta_1(2\tau,5\tau)}.$$
Indeed, since $$\theta_1(v,\tau) = 2 q^{1/4} (\sin \pi 
v)(q^2,q^2)_\infty (q^2e^{-2i\pi v},q^2)_\infty (q^2e^{2i\pi v},q^2)_\infty $$
then $$\theta_1(\tau,5\tau) = i q^{1/4} (q^2,q^{10})_\infty 
(q^8,q^{10})_\infty (q^{10},q^{10})_\infty ;$$ $$ \theta_1(2\tau,5\tau) = i 
q^{-3/4} (q^4,q^{10})_\infty (q^6,q^{10})_\infty (q^{10},q^{10})_\infty .$$
From which we deduce, $$R(q) = q^{-6/10} \frac {(q^2,q^{10})_\infty 
(q^8,q^{10})_\infty}{(q^4,q^{10})_\infty (q^6,q^{10})_\infty} = q^{-6/10} 
\frac {\theta_1(\tau,5\tau)}{\theta_1(2\tau,5\tau)}.$$

The result follows from theorem 2-5 which asserts that
$$\theta_1(v,\tau) = \theta_4(0,\tau)\ \exp[ i\pi (v-\frac {1}{2}+\frac 
{1}{4}\tau) - \sum_{p\geq 1} \sum_{k\geq 0}\frac {1}{p} \bigg( \frac 
{\sin \pi (v+\frac {1}{2}\tau)}{(\sin (k+\frac {1}{2})\pi 
\tau)}\bigg)^{2p}].$$ 
Thus, 
$$R(q) = exp [-i\pi \frac {8\tau}{5} + \sum_{p\geq 1}  \sum_{k\geq 0} - 
\frac {1}{p} \frac {[(\sin \frac {7\pi \tau }{2})^{2p} - (\sin \frac 
{9\pi \tau }{2})^{2p}]}{(\sin (k+\frac {1}{2})5\pi \tau)^{2p}}].$$ 
Another equivalent expression may be obtained by the same techniques
$$ R(q) = exp [-i\pi \frac {8\tau}{5} + \sum_{p\geq 1} c_{2p}(5\tau+1) 
[(\cos \frac {7\pi \tau }{2})^{2p} - (\cos \frac {9\pi \tau 
}{2})^{2p}].$$

\newpage

{\bf References}\\

[A-L] P. Appell, E. Lacour, \quad {\it Fonctions elliptiques et 
applications}\quad Gauthiers-Villard ed., Paris (1922).\\

[A] T.M. Apostol, \quad {\it Modular functions and Dirichlet series in 
number theory}\quad 2nd Ed. New-York; Springer-Verlag, (1997).\\

[B] A.Erdélyi, W.Magnus, F.Oberhettinger,  F.Tricomi, \quad {\it  
Higher transcendental
functions}\quad  Vol. I and III. Based on notes left by H. Bateman. 
Robert E. Krieger Publish. Co.,
Inc., Melbourne, Fla., (1981).\\

[Be] R.E. Belmann, \quad {\it A brief introduction to theta 
functions},\quad New York, Holt, Rinehart and Winston [1961]\\

[C] R. Chouikha \quad \textit{On the expansions of elliptic functions 
and applications}\ in \
Algebraic methods and q-special functions, 
 C.R.M. Proceedings  and Lectures Notes, A.M.S., vol 22 , p.53-57,  
Providence, (1999). \\

[C1] R. Chouikha \quad \textit{Sur des developements de fonctions 
elliptiques} \ Publ. Math., Fac. des Sc. de Besancon. Theorie des nombres, 
vol 8 , p. 1-9,(1989).\\

[C2] R. Chouikha \quad \textit{Note on trigonometric expansions of 
theta functions} \ Proc. of the  6th OPSFA, (2002).\\

[F] R. Fricke\quad {\it Elliptische Funktionen.}\quad  Encycl. der 
Math. Wissensch., t.II, B.3,
Teubner, Leipzig (1924).\\

[L] Z.-G. Liu \quad {\it Two theta functions identities and some 
Eisentein series identities of Ramanujan}\quad  
Rocky Mountain J. of Math., vol 34, 2, p.713-732 (2004).\\

[M] D. Mumford \quad {\it Tata lectures on theta. vol II. Jacobian 
theta functions and differential equations} \quad With  C. Musili, M. Nori, 
E. Previato, M. Stillman and H. Umemura. Progress in Mathematics, 43. 
Birkhäuser, Boston, (1984).\\

[W-W] E.T. Whittaker,G.N. Watson \quad {\it A course of Modern 
Analysis}\\ Cambridge (1963).

\end{document}